\documentclass[a4paper, 12pt,leqno]{article}

\def\publicationcopy{}

\usepackage{ifpdf} 
\ifpdf 
\RequirePackage[colorlinks,hyperindex,plainpages=false,a4paper,%
  pdftitle={On the Representation of Almost Primes by Many Quadratic
  Forms},%
  pdfauthor={Gihan Marasingha},%
]{hyperref}
\else
\RequirePackage[hypertex,a4paper]{hyperref}
\fi

\ifx\publicationcopy\undefined\usepackage[active]{srcltx}\fi

\usepackage{amsmath,amssymb,amsthm,amsfonts,setspace}

\ifx\publicationcopy\undefined\usepackage{color}\fi
\usepackage{verbatim}

\newtheorem{theorem}{Theorem}[section]
\newtheorem{lemma}{Lemma}[section]

\newcommand{\cR}{\mathcal{R}}  

\newcommand{\cU}{\mathcal{U}}
\newcommand{\cA}{\mathcal{A}}
\newcommand{\cRo}{\mathcal{R}^{(0)}}

\newcommand{\fP}{\mathfrak{P}} 
\newcommand{\fPbar}{\overline{\mathfrak{P}}}
\newcommand{\fA}{\mathfrak{A}}

\newcommand{\x}{\mathbf{x}}   
\newcommand{\bd}{\mathbf{d}}
\newcommand{\bc}{\mathbf{c}}
\newcommand{\bdf}{\mathbf{f}}

\newcommand{\bv}{\mathbf{v}}
\newcommand{\by}{\mathbf{y}}
\newcommand{\bz}{\mathbf{z}}
\newcommand{\bQ}{\mathbf{Q}}

\newcommand{\lstard}{\Lambda_{\bd}^{*}}
\newcommand{\lstarc}{\Lambda_{\bc}^{*}}
\newcommand{\ld}{\Lambda_{\bd}}
\newcommand{\rstar}{\rho^{*}}

\newcommand{\R}{\mathbb{R}}  

\newcommand{\N}{\mathbb{N}}
\newcommand{\Z}{\mathbb{Z}}

\newcommand{\eps}{\varepsilon}
\newcommand{\smod}[1]{\;(\text{mod }#1)}  

\newcommand{\ceil}[1]{\lceil{#1}\rceil}
\newcommand{\jacobi}[2]{\genfrac{(}{)}{}{}{#1}{#2}}

\newcounter{saveeqn}  
\newcommand{\alpheqn}{\setcounter{saveeqn}{\value{equation}}
  \setcounter{equation}{0}
  \renewcommand{\theequation} 
     {\mbox{\Alph{equation}}}}
\newcommand{\reseteqn}{\setcounter{equation}{\value{saveeqn}}
  \renewcommand{\theequation}{\arabic{equation}}}

\newenvironment{acknowledgements}{\vspace{1
    ex}\textbf{Acknowledgements.}}{}

\ifx\publicationcopy\undefined\newenvironment{explain}{\color{magenta}}{}\else%
\fi

\DeclareMathOperator{\vol}{vol}
\DeclareMathOperator{\Res}{Res}

\DeclareMathOperator{\lcm}{lcm}

\begin{document}

\title{On the Representation of Almost Primes by Sets of Quadratic
  Forms}

\author{Gihan Marasingha}
\date{20th July, 2006}

\maketitle




\section{Introduction} \label{chap:introduction}

In a previous paper \cite{Mar06}, we investigated an approximation to
Schinzel's Hypothesis H by showing that the product $q_1(\x)
q_2(\x)$ of two irreducible binary quadratic forms has at most 5 prime
factors for infinitely many values $\x$, if the forms satisfy certain
properties.  The present work generalises our ideas to products of
arbitrarily many irreducible binary quadratic forms, and we compute
results which extend the work of Diamond and Halberstam
\cite{DH97}, in which they consider irreducible polynomials.  Indeed,
we have the following Theorem:

\begin{theorem} \label{thm:main}
  Let $q_i(x,y) := a_i x^2 + 2b_i xy + c_i y^2$, for $i=1,\ldots,g$, be
  a set of irreducible quadratic forms over the integers such that
  $a_i \equiv 1 \smod{4}$.  Let
$$
D:= \prod_{p\le 2g} p \prod_{k=1}^g{a_k}{c_k}{\delta_k} \prod_{i<j}
\Res(q_i,q_j),
$$
where $\delta_i$ is the discriminant of the form $q_i$ and
$\Res(q_i,q_j)$ is the resolvent of the forms $q_i$ and $q_j$.  If $D
\ne 0$ and if there exists $\bz \in \Z^2$ such that $(q_i(\bz);D)=1$
for $i=1,\ldots,g$, then
$
q_1(\x) \ldots q_g(\x) = P_{r_M},
$
where $r_M$ is given in Table \ref{table:r_M}, and where $P_n$ denotes a number
with at most $n$ prime factors.  Moreover, if $\cRo$ is a convex
subset of $\R^2$ with piecewise continuously differentiable boundary,
then there exists a positive absolute constant $\gamma_g <1$ such that
for all sufficiently large $X$,
$$
\# \{\x \in X \cRo\, : \,q_1(\x) \ldots q_g(\x) = P_{r_M} \} \gg X^2
\prod_{p<X^{\gamma_g}} \left( 1 - \frac{\omega(p)}{p} \right),
$$
where the implied constant depends at most on the forms and on the
region $\cRo$, and where
$$
\omega(p) = g + \sum_{n=1}^g \jacobi{\delta_i}{p} - \frac{1}{p}\left( g -1 +
  \sum_{n=1}^g \jacobi{\delta_i}{p} \right).
$$
\end{theorem}

\begin{table}[h]
\caption{$r_M$ compared with the values $r_{DH}$ of Diamond
  and Halberstam \cite{DH97}}
  \label{table:r_M}
\begin{tabular}{r| r r r r r r r r r}
$g$ & 2 & 3 & 4 & 5 & 6 & 7 & 8 & 9 & 10 \\
\hline
$r_{DH}$ & 7 & 12 & 17& 23& 29 & 35 & 41 & 47 & 54 \\
\hline
$r_M$ & 5 & 8 & 12 & 16 & 20 & 25 & 29 & 34 & 39
\end{tabular}
\end{table}

Note that Diamond and Halberstam do not consider the problem of almost
primes represented by forms, but one may na\"ively apply their results
on polynomials by fixing one of the variables of the quadratic forms.
The values $r_{DH}$ in Table \ref{table:r_M} give the result of such
an application.

The argument we present will follow the same structure as our previous
paper \cite{Mar06}.  The multidimensional sieve of Diamond and
Halberstam \cite{DH97} will give our values $r_M$, but we have a to
develop a new level of distribution formula to use in the sieve.  In
doing so, we find a new type of upper bound for the quantity
$\rstar$.  The verification of sieve condition
(\ref{eqn:rkappaalpha}) presents a departure 
from the original proof, as we find a way to construct an appropriate
set of mutually pairwise coprime integers from a more general set of
squarefree numbers.
\subsection{Notation}

$(a_1;\ldots;a_n)$ is the highest common factor of the numbers $a_i$,
for
$i=1,\ldots,n$; $P_n$ represents a number with at most $n$ prime
factors; $d(n)$ is
the number of positive divisors of $n$; $\phi(n)$ is the number of
nonnegative integers less than and prime to $n$; $\mu(n)$ is the
M\"obius function; $\nu(n)$ is the number of distinct prime factors of
$n$; $\jacobi{\cdot}{p}$ is the Legendre symbol; and $\delta_{ij}$ is the
Kronecker delta function (that is $\delta_{ij}=0$ if $i \ne j$ and
$\delta_{ij} = 1$ if $i=j$).  We will use the symbol $C$ to denote a
positive numerical constant, though its value may vary in the course
of a proof.


\section{The Level of Distribution} \label{chap:lod}

With $q_i$, $\bz$, and $D$ as in the statement of Theorem
\ref{thm:main}, we make the following definitions:
\begin{align*} 
\ld &:= \{ \x \in \Z^2\, : \,d_i | q_i(\x) \quad(i=1,\ldots,g) \},&
\rho(\bd) &:= \# (\ld \cap [0,d_1\ldots d_g)^2), \\
\lstard  &:= \{ \x \in \ld\, : \, (\x;d_1\ldots d_g)=1\}, &
\rstar(\bd) &:= \# (\lstard \cap [0,d_1\ldots d_g)^2), \\
\Psi &:= \{ \x \in \Z^2\, :\, \x \equiv \bz \pmod{D} \}, & \Psi_b &:=
\{ \x \in \Z^2\, : \, b \x \in \Psi\}.
\end{align*}

\begin{theorem}[Level of Distribution] \label{thm:lod}
Let $q_i(\x)$ and $D$ be defined as in Theorem \ref{thm:main}.  For
any real numbers $M,Q_1,\ldots,Q_g\ge 1$, let
$$
T(M,\bQ) := \sum_{\substack{d_i \le Q_i \\ (d_i;D)=1}}
\sup_{\partial(\cR) \le  M} \left| \# (\ld \cap \cR \cap \Psi) - \frac{\vol(\cR)
    \rho(d_1,\ldots,d_g)}{(d_1 \ldots d_g D)^2} \right|.
$$
Then, writing $Q := Q_1\ldots Q_g$, there exist constants
$\nu_1, \nu_2 >1 $, depending only on $g$  such that
$$
T(M,\bQ) \ll Q (\log 2 Q)^{\nu_1} + M \sqrt{Q} (\log 2 Q)^{\nu_2}.
$$
\end{theorem}

\subsection{Transition from $\lstard$ to $\ld$} \label{sec:transition}

We begin with the following bridging result, which will be employed in
Section \ref{sec:lodunstarred} to
express the unstarred sum in
terms of the starred sum, leading to Theorem
\ref{thm:lod}.

\begin{lemma}[Transition Formula] \label{lem:transition}
Let $D \in \N$ and suppose that $(d_i;D)=1$ for $i=1,\ldots,g$.  Then
we have
$$
  \# \left( \ld \cap \cR  \cap \Psi \right) = \sum_{b | \psi(\bd)} \#
  \left( \lstarc \cap \cR / b  \cap \Psi_b\right),
$$
where $c_i := d_i / (d_i ; b^2)$, for $i=1,\ldots,g$, and the
multiplicative function $\psi$ is defined by
$
\psi(p^{\alpha_1},\ldots,p^{\alpha_g}) :=
p^{\ceil{\max(\alpha_1,\ldots,\alpha_g)/2}}.
$
\end{lemma}

By definition, $\# (\ld \cap \cR \cap \Psi) = \# \{ \x \in \cR : d_i | q_i(\x),\
 i=1,\ldots,g,\ \x \in \Psi \}$. As in our pairs of forms article
 \cite{Mar06}, the Lemma follows by partitioning
 this set according to $(\x;\psi(\bd))$.
\begin{explain}
Such a partition lets us write  $\# (\ld \cap \cR \cap \Psi)$ as:
$$
 \sum_{b | \psi(\bd)} \# \{ \x \in \cR : q_i(\x) \equiv 0
 \smod{d_i},\ (\x ; \psi(\bd)) = b, \ \x \in \Psi \}.
$$
We claim that if $b | \psi(\bd)$ then $\psi(\bd)/b =
\psi(\bc)$. Indeed, if $b = \prod p^\beta$, and $d_i = \prod
p^{\alpha_i}$, then $\psi(\bd)/b = 
\prod_p p^{\lceil (\max(\alpha_1,\ldots,\alpha_g) -2 \beta)/2\rceil }$.
On the other hand, $$\psi(\bc) = \prod_p p^{\lceil \max_{i=1,\ldots,g}(\alpha_i -
  \min(\alpha_i,2\beta))/2 \rceil}.$$
Now for each $i$, $$\alpha_i - \min(\alpha_i,2\beta) = \begin{cases}0 &
  \text{ if } \alpha_i \le 2 \beta, \\ \alpha_i - 2 \beta & \text{ if }
  \alpha_i \ge 2 \beta.\end{cases}$$
Therefore $$\max_{i=1,\ldots,g}(\alpha_i - \min(\alpha_i,2\beta)) =
\max_{i=1,\ldots,g}(\alpha_i)-2\beta.$$

It is clear that $(\by;\psi(\bc))=1 \iff (\by;c_1\ldots c_g)=1$.
Write $\x = b \by$; we have $b | \psi(\bd)$, so $(b;D)=1$.  We may
write $b^{-1} \x \equiv \by \smod{D}$, and deduce
$$
\# (\ld \cap \cR \cap \Psi)  = \sum_{b|\psi(\bd)} \# \{ \by \in \cR/b
: c_i |q_i(\by), \ (\by;c_1\ldots c_g)=1, \by \in \Psi_b\},
$$
hence the result.
\end{explain}
\subsection{Upper Bounds for $\rho$}

\subsubsection{The Function $\rstar$}

\begin{lemma} \label{lem:rstarbound}
  For every prime $p$ one has that $\rstar(p^{e_1},\ldots, p^{e_g})
  \ll p^{\max(e_1,\ldots,e_g)},
 $
and if $(p;D)=1$ then
 \begin{equation} \label{eqn:rstarsplit}
   \rstar(p^{e_1},\ldots,p^{e_g}) = 0, \text{ if } e_i, e_j>0 \text{
     for some }i \ne j.
 \end{equation}
\end{lemma}

Assume $\rstar(p^{e_1},\ldots,p^{e_g}) \ne 0$.  As in \cite{Mar06}, we
may demonstrate that if $i \ne j$, then
$p^{\min(e_i,e_j)}|\Res(q_i,q_j)$, and hence $p^{\min(e_i,e_j)}|D$.
\begin{explain}
  To see this, note that there exists $\x
= (x_1,x_2)$ such that $(\x;p)=1$ and $p^{e_i} | q_i(\x)$ for
$i=1,\ldots g$.  Without loss of generality, $p \nmid x_2$.  Define
$Q_i(Y) := q_i(Y,1)$, and $y \equiv x_1 x_2^{-1}
\smod{p^{\max(e_1,\ldots,e_g)}}$.  Then $0 \equiv q_i(x_1,x_2) \equiv
x_2^2 Q_i(y) \smod{p^{e_i}}$.  Therefore $p^{e_i} | Q_i(y)$ for
$i=1,\ldots, g$.  So for all $i \ne j$, $Q_i(\by) \equiv Q_j(\by)
\equiv 0 \smod{p^{\min(e_i, e_j)}}$.  We deduce
$p^{\min(e_i,e_j)} | \Res(Q_i,Q_j)=\Res(q_i,q_j)$.
\end{explain}
Suppose $(p;D)=1$; then $\min(e_i,e_j)=0$ for every pair $i\ne j$.
This suffices for the second part of the Lemma.  The first part is
proved by induction on $g$.  Our previous paper \cite{Mar06} contains
a proof for the base case, $g=2$.  In general, one has
\begin{align*}
\rstar(p^{e_1},\ldots,p^{e_g}) &\le \# \{ \x \smod{p^{e_1+\ldots + e_g}}
: p^{e_i}| q_i(\x), (\x;p)=1, i=2,\ldots,g\} \\
& = p^{2e_1} \# \{ \x \smod{p^{e_2+\ldots+e_g}} : p^{e_i} | q_i(\x),
 i=2,\ldots, g ; (\x;p)=1\} \\
& \ll p^{2 e_1} p^{\max(e_2,\ldots,e_g)} \text{ by induction on }g \\
& \le p^{2 e_1 \max(e_1,\ldots,e_g)} 
\end{align*}
Generalising, we have that for all $i$, $\rstar(p^{e_1},\ldots,
p^{e_g}) \ll p^{2e_i} p^{\max(e_1,\ldots,e_g)}$.  Thus
$$
\rstar(p^{e_1},\ldots,p^{e_g}) \ll p^{2 \min(e_i, e_j)}
p^{\max(e_1,\ldots e_g)} \ll D^2 p^{\max(e_1,\ldots,e_g)} \ll p^{\max(e_1,\ldots,e_g)}.
$$

\subsubsection{The Function $\rho$}

\begin{lemma} \label{lem:rbound}
  Let $p$ be a prime and let $e_1, \ldots, e_g$ be non-negative
  integers.  Let $\sigma$ be a permutation in $S_g$ such that
  $e_{\sigma(1)} \le \ldots \le e_{\sigma(g)}$.  Then
$$
\rho(p^{e_1},\ldots,p^{e_g}) \ll (e_{\sigma(g)} - e_{\sigma(g-1)} +1 )p^{2
e_{\sigma(1)}+\ldots+2e_{\sigma(g-1)}+e_{\sigma(g)}}.
$$
Also,
$$\rho(p,\ldots,p) = \rstar(p,\ldots,p)+p^{2(g-1)} = p^{2(g-1)} + O(p),$$
and for all but a finite set of primes $p$, one has
$$
\rho(p,1,\ldots,1) \le 2p ;\ \rho(1,p,1,\ldots,1) \le 2p,\ \ldots, \
\rho(1,\ldots,1,p) \le 2p.$$
\end{lemma}

An application of Lemma \ref{lem:transition} provides us with the
following formul\ae:
\begin{equation}
  \rho(\bd) = \sum_{b | \psi(\bd)} \# (\lstarc \cap [0,d_1 \ldots d_g /
  b)^2) = \sum_{b | \psi(\bd)} \rstar(\bc) \left( \frac{(d_1 ; b^2)
      \ldots (d_g;
  b^2)}{b} \right)^2, \label{eqn:rhotransition}
\end{equation}
where $c_i = d_i / (d_i ; b^2)$.
In particular, $\rho(p^{e_1},\ldots,p^{e_g})$ can be written as
\begin{equation}
  \sum_{0 \le \beta \le \ceil{e_g/2}} \rstar \left(
  \frac{p^{e_1}}{(p^{e_1} ; p^{2 \beta})}, \ldots,
  \frac{p^{e_g}}{(p^{e_g} ; p^{2 \beta})}
  \right) \left( \frac{(p^{e_1} ; p^{2 \beta})\ldots(p^{e_g} ; p^{2
  \beta})}{p^\beta} \right)^2. \label{eqn:rhorstar}
\end{equation}

Using these formul\ae and our upper bound for $\rstar$,
we may derive a proof of the first two statements, much as in \cite{Mar06}.
\begin{explain}
Indeed, without loss of generality, we shall assume $e_1 \le \ldots \le e_g$.
Split the range of summation as $0 \le 2 \beta < e_1$, $e_1 \le 2 \beta <
e_2, \ldots, e_{g-1} \le 2 \beta < e_g$ and $\beta=\ceil{e_g/2}$.  We
have the following expression for
$\rho(p^{e_1},\ldots, p^{e_g})$:
\begin{align*}
&\sum_{0 \le \beta < e_1/2}
  \rstar(p^{e_1-2\beta},\ldots,p^{e_g-2\beta}) p^{\beta(4g-2)} \\
&+ \sum_{e_1/2 \le \beta < e_2/2}
 \rstar(1,p^{e_2-2\beta},\ldots,p^{e_g-2\beta}) p^{2 e_1}
 p^{\beta(4g-6)} \\ 
&+ \sum_{e_2/2 \le \beta < e_3/2}
 \rstar(1,1,p^{e_3-2\beta},\ldots,p^{e_g-2\beta}) p^{2e_1 +2e_2}
 p^{\beta(4g-10)} \\
&+\ldots+ \sum_{e_{g-1}/2 \le \beta < e_g/2}
   \rstar(1,\ldots,1,p^{e_g-2\beta}) p^{2e_1+\ldots 2e_{g-1}}
   p^{2\beta} \\
&+ p^{2 e_1+\ldots + 2e_{g-1}+e_g}.
\end{align*}
Via an application of our upper bound for $\rstar$, we see that for
each $0\le i \le g-2$, the 
sum over the range $e_i/2 \le \beta < e_{i+1}/2$ is bounded by
$p^{2e_1+\ldots+2e_i+e_g}p^{2e_{i+1}(g-i-1)}$.  This is an increasing
function of $i$, so $\rstar$ is bounded above by $p^{2e_1+\ldots+2
  e_{g-1}+e_g} + (e_g-e_{g-1})p^{2e_1+\ldots+2e_{g-1}+e_g}$, as
required.

An application of equation (\ref{eqn:rhorstar}) gives
$$
\rho(p,\ldots,p) = \rstar(p,\ldots,p)+p^{2(g-1)} = p^{2(g-1)} +O(p),
$$
\end{explain}
The final statement of Lemma \ref{lem:rbound} is proved using the
one-form result of \cite{Mar06}.

\subsection{Level of Distribution---Starred Version}

\begin{lemma} \label{thm:lodstar}
Define
  \[
  T^{*}(M,\bQ) := \sum_{\substack{d_i \le Q_i \\ (d_i;D)=1}} \sup_{\cR
  : \partial R \le M}
  \left| \# (\lstard \cap \cR \cap \Psi) - \frac{\rstar(d_1,\ldots, d_g)}
  {(d_1 \ldots d_g D)^2}
  \vol(\cR)\right|.
\]
Then there exist constants $\nu_1'$
and $\nu_2'$ depending only on $g$ such that
\[ T^{*}(M,\bQ) \ll M \sqrt{Q} ( \log(2 Q))^{\nu_1'} +
Q ( \log (2 Q) )^{\nu_2'}
\]
uniformly for $M > 0$ and $Q_1, \ldots, Q_g \ge 1$.
\end{lemma}

\subsubsection{The Quantities $\lstard$}

Assume that $\bd = (d_1, \ldots, d_g)$ is fixed and
define $a := d_1 \ldots d_g$.  Let $\cU(a)$ be the set of equivalence classes
of $\x \in \Z^2$ under multiplication with $(x_1 ; x_2 ; a)=1$, as in
\cite{Mar06}.  For a given $\cA \in \cU(a)$, one has that $\cA \subset
\lstard$ or $\cA \cap \lstard = \emptyset$.  Define $\cU'(\bd) := \{
\cA \in \cU(d_1\ldots d_g) : \cA \subset \lstard \}.$  Then $\lstard$
is the disjoint union
$
\lstard = \bigcup_{\cA \in \cU'(\bd)}{\cA}.
$

Now $\# (\cA \cap [0,a)^2) = \phi(a)$, so
\begin{equation} \label{eqn:rhoeqn}
 \rstar(\bd) = \# \cU'(\bd) \phi(d_1 \ldots d_g).
\end{equation}

Using this equation, the summand can be bounded from above:
\begin{equation} \label{eqn:lodstardd}
\begin{split}
  & \left| \# (\lstard \cap \cR \cap \Psi) - \frac{\rstar(\bd)}
  {(d_1 \ldots d_g D)^2}
  \vol(\cR)\right|  \\
   \le & \sum_{ \cA \in \cU'(\bd)} \left| \# (\cA \cap \cR \cap \Psi) -
  \frac{\phi(d_1 \ldots d_g)}{(d_1 \ldots d_g D)^2} \vol(\cR) \right|.
\end{split}
\end{equation}

\subsubsection{Estimating $\#(\cA \cap \cR \cap \Psi)$}

Choose
$\cA \in \cU(a)$ and define $G(\cA)$ by:
\[
G(\cA) := \{ \x \in \Z^2 : (\exists \lambda \in \Z) (\exists \by \in
\cA) (\x \equiv \lambda \by \smod{a}) \}. 
\]
As in \cite{Mar06}, the lattice $G(\cA)$ has a non-zero element of minimal
length $\bv(\cA)$.  Moreover, $\bv(\cA)$ satisfies $|\bv(\cA)|^2 \ll
a$, and
$$
\#(\cA \cap \cR \cap \Psi)
= \frac{\phi(a)}{a^2D^2} \vol(\cR) + O \left( d(a) \left(
     \frac{M}{|\bv(\cA) |} + 1 \right) \right). \label{eqn:caeqn}
$$
Substituting this into equation (\ref{eqn:lodstardd}), we
have:
\begin{equation} \label{eqn:lodstar5}
  \begin{split}
      T^{*}(M, \bQ) & \ll \sum_{\substack{d_i \le Q_i \\ (d_i;D)=1}}
      \sum_{\cA \in
      \cU'(\bd)} d( d_1 \ldots d_g) \left(
      \frac{M}{|\bv(\cA)|} + 1 \right) \\
     & \ll M \sum_{d_i \le Q_i}  d(d_1 \ldots d_g) \sum_{\cA \in
       \cU'(\bd)}
       \frac{1}{| \bv(\cA) |} + \sum_{d_i \le Q_i} d(d_1 \ldots d_g)       
       \# \cU'(\bd) \\
     & = M T_1^{*}(\bQ) + T_2^{*}(\bQ), \text{ say.}
  \end{split}
\end{equation}

\subsubsection{Evaluating \protect{$T_1^{*}(\bQ)$}}

We shall prove:
\begin{lemma} \label{lem:t1bound}
The quantity
$$
T_1^*(\bQ) := \sum_{d_i \le Q_i} d(d_1 \ldots d_g) \sum_{\cA \in \cU'(\bd)}
       \frac{1}{| \bv(\cA) |}
$$
satisfies the upper bound $ T_1^*(\bQ) \ll \sqrt{Q} ( \log 2 Q)^{\nu_1'},$
for some constant $\nu_1'$ depending only on $g$, and where $Q := Q_1
\ldots Q_g$.
\end{lemma}

For $\cA \in \cU'(\bd)$, we have $| \bv(\cA)| \ll \sqrt{d_1 \ldots
  d_g}$.  Consequently,
$$
  T_1^{*}(\bQ) \le \sum_{0 < | \bv| \ll \sqrt{Q}} \frac{1}{| \bv
  |} \sum_{\substack{d_i \le Q_i \\ d_i | q_i(\bv)}} d(d_1 \ldots d_g)
\# \cU'(\bd).
$$

Write $d_i = \prod p^{e_i}$ and apply equation (\ref{eqn:rhoeqn}):
$$
  \# \cU'(\bd) = \frac{\rstar(\bd)}{\phi(d_1 \ldots d_g)}
  = \prod_p  \frac{\rstar(p^{e_1},\ldots,p^{e_g})}{\phi(p^{e_1+\ldots+e_g})},
$$
by multiplicativity of $\rstar$.  By  Lemma \ref{lem:rstarbound}, if
$(p;D)=1$ and if at least two of the 
$e_i$ are positive,
then $\rstar(p^{e_1},\ldots,p^{e_g})=0$.  Thus $\# \cU'(\bd)=0$ unless
for all $p$
 satisfying $(p;D)=1$, we have $e_i=0$ for all but at most one $i$.  In which
case, we may write:
$$
  \# \cU' (\bd)  \le
  \prod_{\substack{p : \\(p;D)=1}}
  \prod_{i=1,\ldots,g}
  \frac{\rstar(p^{e_i \delta_{1i}},\ldots,p^{e_g \delta_{g_i}})}
       {\phi(p^{e_1 \delta_{1i}+\ldots e_g \delta_{gi}})}
  \prod_{\substack{p  : \\ p | D}} C
   \frac{p^{\max(e_1,\ldots,e_g)}}{\phi(p^{e_1+\ldots e_g})} \ll 2^{\nu(d_1\ldots d_g)},
$$
much as in \cite{Mar06}.

\begin{explain}
Observing that $2^{\nu(d_1 \ldots d_g)} \le d(d_1 \ldots d_g)$, and using the
submultiplicativity property of $d$, we deduce that
$
d(d_1 \ldots d_g) \# \cU'(\bd) \ll d(d_1)^2 \ldots d(d_g)^2,
$
whence
$$
T_1^{*}(\bQ) \ll \sum_{0 < | \bv| \le \sqrt{Q}} \frac{1}{|\bv|}
  \sum_{d_i | q_i(\bv)} d(d_1)^2 \ldots d(d_g)^2.
$$
\end{explain}
Defining the function $h$ by
$
h(n) := \sum_{a | n} d(a)^2,
$
we have
\begin{equation} \label{eqn:lodstar6}
  T_1^{*}(\bQ) \ll \sum_{0 < |\bv| \le \sqrt{Q}} \frac{1}{|\bv|}
  h(q_1(\bv)) \ldots h(q_g(\bv)).
\end{equation}
Regarding this function, we have the following Lemma, to be found
in \cite{Mar06}:

\begin{lemma} \label{lem:hdivisor}
The function $h$ is multiplicative. Moreover, $h$ is
submultiplicative in the sense that for all $m_1, m_2$, we have
$h(m_1 m_2) \le h(m_1) h(m_2)$.  Furthermore, we have $h(p) \ll 1$
uniformly in $p$.
Let $\eta \ge 1$, then for every natural number $n$, there exists a
positive integer $m$ satisfying $m|n$, $m \le n^{1 /\eta}$, and
$$
 h(n) \ll_\eta h(m)^{1+\lfloor \eta\rfloor}.
$$
\end{lemma}

We apply Lemma \ref{lem:hdivisor} to equation
(\ref{eqn:lodstar6}), with $\eta = 2g$, to obtain:
$$
T_1^*(\bQ) \ll \sum_{\substack{j \ge 0 \\ P = 2^j \ll \sqrt{Q}}}
\frac{1}{P} \sum_{P \le |\bv| \le 2P} \sum_{\substack{m_i | q_i(\bv)
    \\ m_i \ll P^{1/g} \\ i=1,\ldots,g}} h(m_1)^{2g+1} \ldots h(m_g)^{2g+1};
$$
\begin{explain}
where we have split the range for $| \bv |$ into dyadic intervals, and
used the fact that $q_i(\bv) \ll | \bv|^2$ to deduce that $m_i \le
|q_i(\bv)|^{1/(2g)}$ implies $m_i \ll P^{1/g}$;
\end{explain}
we have
$$
T_1^*(\bQ) \ll \sum_{\substack{j \ge 0 \\ P = 2^j \ll \sqrt{Q}}}
\frac{1}{P}
\sum_{m_i \ll P^{1/g}} h(m_1)^{2g+1} \ldots h(m_g)^{2g+1} \sum_{\substack{|\bv| \le
    2P \\ q_i(\bv) \equiv 0 \smod{m_i} \\ i=1, \ldots g}} 1.
$$
The innermost sum
\begin{explain}
is of order $\rho(m_1,\ldots, m_g)\{ P^2 / (m_1
\ldots m_g)^2 +
P/(m_1 \ldots m_g)\}$,
where the second term accounts for the error at the boundary.  We have
arranged that $m_1 \ldots m_g \ll P$, whence the second term is subsumed by
the first, and hence the innermost sum
\end{explain}
is of order $\rho(m_1, \ldots, m_g) P^2
/ (m_1 \ldots m_g)^2$, so
\begin{align*}
T_1^*(\bQ) & \ll \sum_{\substack{j \ge 0 \\ P=2^j \ll \sqrt{Q}}} P
\sum_{m_i \ll P^{1/g}}\frac{h(m_1)^{2g+1} \ldots h(m_g)^{2g+1}
  \rho(m_1, \ldots, m_g)}{m_1^2 \ldots m_g^2}\\ 
 & \ll \sqrt{Q} \sum_{m_i \ll Q^{1/(2g)}}
 \frac{h(m_1)^{2g+1} \ldots h(m_g)^{2g+1} \rho(m_1, \ldots,
   m_g)}{m_1^2 \ldots m_g^2}.
\end{align*}
The summand is multiplicative, so we have the following upper bound:
$$
T_1^*(\bQ) \ll \sqrt{Q} \prod_{p \ll Q^{1/(2g)}}
\sum_{\substack{e_i=0 \\i=1,\ldots,g}}^\infty \frac{h(p^{e_1})^{2g+1}
  \ldots h(p^{e_g})^{2g+1} \rho(p^{e_1},\ldots,p^{e_g})}{p^{2e_1+\ldots+2e_g}}.
$$
Let $k(e_1,\ldots,e_g)$ denote the summand.  Using our upper bound for
$\rho$, we estimate the
sum $S= \sum_{e_i=0}^\infty k(e_1,\ldots,e_g)$ as follows:
\begin{align*}
S \le & \sum_{e_1=0}^\infty \sum_{e_2,\ldots,e_g \le e_1}
k(e_1,\ldots,e_g) + \sum_{e_2=1}^\infty \sum_{\substack{e_i \le e_2 \\
  i \ne 2}} k(e_1,\ldots,e_g) \\ & + \ldots + \sum_{e_g=1}^\infty
\sum_{\substack{e_i \le e_g \\ i \ne g}} k(e_1,\ldots,e_g).
\end{align*}

Now
$$
\sum_{e_1=0}^\infty \sum_{\substack{e_i\le e_1 \\ i \ne 1}}
k(e_1,\ldots,e_g) = 1 + \frac{C_1}{p}+\frac{C}{p^2} +
\sum_{e_1=2}^\infty \sum_{\substack{e_i \le e_1 \\ i \ne 1}}
k(e_1,\ldots,e_g).
$$
\begin{explain}
The first term in this expression arises from $k(0,\ldots,0)$, the
second from $k(1,0,\ldots,0)$, and the third from those sets of $e_i$
where $e_1=1$ and at least one other $e_i=1$, the bound deriving from
the $n$-variable result $\rho(p,\ldots,p) \ll p^{2(n-1)}$.
\end{explain}
We have
$$
S' := \sum_{e_1=2}^\infty \sum_{\substack{e_i \le e_1 \\ i \ne 1}}
k(e_1,\ldots,e_g) \ll \sum_{e_1=2}^\infty \sum_{e_i \le e_1}
\frac{e_1^{6g+3} \ldots e_g^{6g+3} e_1}{p^{e_1}},
$$
using $h(p^e) \ll e^3$ and $\rho(p^{e_1},\ldots,p^{e_g}) \ll e_1 p^{2
  e_2 + \ldots 2 e_g +e_1}$.
So $S' \ll \sum_{e_1=2}^\infty \frac{e_1^A}{p^{e_1}}$, for some $A \in
\N$, depending on $g$.
Therefore $p^2 S' \ll \sum_{e_1=2}^\infty \frac{e_1^A}{p^{e_1-2}} \le
\sum_{e_1=2}^\infty \frac{e_1^A}{2^{e_1-2}} < \infty,$ by the ratio
test.

Thus
$$
\sum_{e_1=0}^\infty \sum_{\substack{e_i\le e_1 \\ i \ne 1}}
k(e_1,\ldots, e_g) \le 1 + \frac{C_1}{p}+\frac{C}{p^2},
$$
for a new constant $C$.  By a similar argument,
$$
\sum_{e_i=1}^\infty \sum_{\substack{e_j \le e_i \\ j \ne i}}
k(e_1,\ldots,e_g) \le \frac{C_1}{p} + \frac{C}{p^2},
$$
for $i=1,\ldots,g$.  Thus, with a possible change of constants $C_1$
and $C$,
$$
\sum_{\substack{e_i=0 \\ i=1,\ldots,g}}^\infty k(e_1,\ldots,e_g) \le 1 +
\frac{C_1}{p}+\frac{C}{p^2}. 
$$

The constant $C_1$ could potentially depend on the forms in question.
However, by a more careful analysis, one can
remove this dependence and ensure that $C_1$ depends only on $g$.

We have:
$$
  T_1^*(\bQ) \ll \sqrt{Q} \prod_{p \ll Q^{1/(2g)}} \left( 1
    + \frac{C'}{p} + \frac{C}{p^2} \right ).  
$$
Much as in \cite{Mar06}, this implies our desired upper bound.
\begin{explain}
Indeed consider the following Lemma:
\begin{lemma}  \label{lem:mertensresult}
  Let $Q > 1$ and $C > 0$ be real numbers.  Let $k$ be a natural
  number and define
$$ S' = \prod_{p \le Q} \left( 1 + \frac{k}{p} +
  \frac{C}{p^2} \right).$$
Then
$$
S' \ll_{k, C}(\log Q)^k.
$$
\end{lemma}
\end{explain}

\subsubsection{Evaluating \protect{$T_2^{*}(\bQ)$}}
We shall prove
\begin{lemma} \label{lem:t2bound}
  The quantity $T_2^*$ satisfies the upper bound
$$
T_2^*(\bQ) \ll Q \log(2 Q)^{\nu_2'}.
$$
for some constant $\nu_2'$, depending only on $g$.
\end{lemma}

\begin{explain}
Recall
$$
T_2^*(\bQ) := \sum_{\substack{d_i \le Q_i \\ i = 1,\ldots,g}}d(d_1
\ldots d_g) \# \cU'(\bd).
$$
\end{explain}
In our analysis of the sum $T_1^*(\bQ)$, we demonstrated that $\#
\cU'(\bd) \ll 2^{\nu(d_1 \ldots d_g)}$.  We have that $2^{\nu(a)} \le
d(a)$ for any $a$ and that the $d$ function satisfies $d(ab) \le d(a)
d(b)$ for any $a$ and $b$. Thus,
\begin{align*}
T_2^*(\bQ) & \ll \left( \sum_{d_1 \le Q_1} d(d_1)^2 \right) \ldots \left(
  \sum_{d_g \le Q_g}d(d_g)^2\right) \ll Q (\log Q_1)^3 \ldots (\log
Q_g)^3 \\
& \ll Q (\log 2 Q)^{3g},
\end{align*}
where we use the AM--GM inequality in the last line.
This proves Lemma
\ref{lem:t2bound}.  Combining this with Lemma \ref{lem:t1bound}
gives us our starred level of distribution formula, Lemma
\ref{thm:lodstar}.

\subsection{Level of Distribution---Unstarred Version}
\label{sec:lodunstarred}

Recall our convention that the symbol $c_i$ represents $d_i / (d_i ;
b^2)$.  We apply Lemma \ref{lem:transition} and equation
(\ref{eqn:rhotransition}) to give the following expression for $T(M,\bQ)$:
\begin{align*}
&\sum_{\substack{d_i \le Q_i \\ (d_i ; D)=1}}
\sup_{\substack{\cR: \\ \partial(\cR) \le M}}
\left| \sum_{b | \psi(\bd)} \left\{ \#(\lstarc \cap \cR / b \cap \Psi_b) -
  \frac{\rstar(\bc)}{(c_1 \ldots c_g D)^2} \vol(\cR/b)  \right\} \right| \\
 &\le \sum_{\substack{c_i \le Q_i  \\ (c_i;D)=1 }} \sum_{\substack{b
     \le Q \\ (b_i ; D ) =1}} \delta(\bQ,\bc,b)
  \sup_{\substack{\cR : \\ \partial(\cR) \le M}}
   L(\bc,b,\cR),
\end{align*}
where we write $L(\bc,b,\cR)$ in place of
$$
\left| \#(\lstarc \cap \cR / b \cap \Psi_b) -
  \frac{\rstar(\bc)}{(c_1 \ldots c_g D)^2} \vol(\cR/b) \right|,
$$
and where $\delta(\bQ,\bc,b) = \# \{ (d_1, \ldots d_g) : d_i \le Q_i,\ c_i = d_i /
(d_i ; b^2),\ b | \psi (\bd) \}$.

We shall derive an upper bound for $\delta$.  Note that
$\delta(\bQ,\bc,b) \le \prod_{i=1}^g \Delta(c_i,b)$, where $\Delta(c_i,b) := \#
\{d_i: c_i = d_i/(d_i;b^2) \}.$  Suppose $d_i = \prod p^{\alpha_i}$,
$b = \prod p^\beta$, and $c_i = \prod p^{\gamma_i}$.  We can estimate
$\Delta(c_i,b)$ by determining how many choices there are for each
$\alpha_i$.  We require $\alpha_i - \min(\alpha_i, 2 \beta) =
\gamma_i$.  For each $i$, we have either $\gamma_i=0$ or
$\gamma_i>0$.  If $\gamma_i=0$, then $\alpha_i \le 2 \beta$, so there
are $2 \beta+1$ choices for $\alpha_i$.  If $\gamma_i >0$, then
$\alpha_i = \gamma_i+2\beta$, so there is only one choice for
$\alpha_i$.  In either case, there are at most $2\beta+1$ choices for
$\alpha_i$.  Thus $\Delta(c_i,b) \le \prod_{p^\beta || b }3 \beta$,
and hence $\delta(\bQ,\bc,b) \le \prod_{p^\beta||b} (3 \beta)^g
=: \theta(b)$.

By analogy with \cite{Mar06}, we have $\sum_{b \le B}
\delta(\bQ,\bc,b) \ll B (\log B)^{A_1}$, where $A_1 := 3^g-1$.
\begin{explain}
Indeed, by
induction on $\beta$, one has $3 \beta \le \binom{\beta+2}{2}$ (for
the induction step, note $\binom{\beta+3}{2} = \binom{\beta+2}{2} +
\binom{\beta+2}{1}).$
One has $d_3(b) = \prod_{p^\beta || b} \binom{\beta+2}{2}$  (see \cite{IK04},
section 1.4) and $\sum_{b \le B} d_3(b)^g \ll B (\log B)^{3^g -1}$
(ibid. section 1.6), so $\sum_{b \le B} \delta(\bQ,\bc,b) \le \sum_{b
  \le B} \theta(b) \ll B (\log B)^{A_1}$.
\end{explain}
If $\delta(\bQ,\bc,b) \ne 0$, then we may deduce that $c_1 \ldots c_g
b \le Q$, restricting the range of summation in the level of
distribution formula.
\begin{explain}
Indeed, suppose that $\delta(\bQ,\bc,b)
\ne 0$, then there exist $d_1, \ldots d_g$ such that $b | \psi(\bd)$ and
$c_i = d_i / (d_i ; b^2)$.  One may verify, by restricting to prime
powers, that $b | \psi(\bd)$ implies $b| (d_1 ; b^2)\ldots (d_g ;
b^2)$.  This may be rewritten as $c_1 \ldots c_g b | d_1 \ldots d_g$,
from which it follows that $c_1 \ldots c_g b \le Q$.
\end{explain}

Our sum $T(M,\bQ)$ is estimated by
\begin{explain}
\begin{align*}
T(M,\bQ) & \ll \sum_{\substack{c_i: \\ (c_i;D)=1 \\ c_i \le Q_i}}
    \sum_{\substack{b : \\ (b;D)=1 \\ c_1 \ldots c_g b \le
        Q}}\delta(\bQ, \bc, b) \sup_{\substack{\cR : \\\partial(\cR)
        \le M }}  L(\bc,b,\cR) \\
& \le \sum_{\substack{j_i : \\ C_i = 2^{j_i} \le Q_i}}
    \sum_{\substack{C_i \le c_i \le 2 C_i \\ (c_i;D)=1}}
    \sum_{\substack{b : \\ (b;D)=1 \\b \le
    \frac{Q}{c_1 \ldots c_g}}} \theta (b)
    \sup_{\substack{\cR : \\ \partial(\cR) \le M }}
   L(\bc,b,\cR).
\end{align*}

If we further split the range for $b$ into dyadic intervals, then
\end{explain}
$$
T(M,\bQ) \ll \sum_{\substack{j_i : \\ C_i = 2^{j_i} \le Q_i}}
    \sum_{\substack{C_i \le c_i \le 2 C_i \\ (c_i;D)=1}}
    \sum_{\substack{k: \\ B = 2^k \le \frac{Q}{c_1 \ldots c_g}}}
    \sum_{\substack{b : \\ (b;D)=1 \\B \le b \le 2 B}} \theta(b)
    \sup_{\partial(\cR) \le M} L(\bc,b,\cR).
$$
\begin{explain}
Our aim is to use the estimate for $\sum_b \theta(b)$, but we need to
handle sensitively the factor of $\sup L(\bc,b,\cR)$.
\end{explain}
For each choice
of $B$, define $b(B)$ by requiring $B \le b(B) \le 2B$, $(b(B);D)=1$
and requiring that for all $b$ with $B \le b \le 2B$ and $(b;D)=1$, one has
$$
\sup_{\partial(\cR) \le M} L(\bc,b,\cR) \le \sup_{\partial(\cR) \le M
} L(\bc,b(B),\cR).
$$
Let $S$ denote the set of integers $B$
such that there are no $b$ in the range $B \le b \le 2B$ with
$(b;D)=1$.  We have the upper bound:
\begin{align*}
T(M,\bQ)
& \ll \sum_{\substack{j_i : \\ C_i = 2^{j_i} \le Q_i}}
    \sum_{\substack{C_i \le c_i \le 2 C_i \\ (c_i;D)=1}}
    \sum_{\substack{k: \\ B = 2^k \le \frac{Q}{c_1 \ldots c_g} \\ B
        \not\in S}}
   B (\log 2B)^{A_1}  \sup_{\partial(\cR) \le M} L(\bc,b(B),\cR) \\
& \le \sum_{\substack{j_i : \\ C_i = 2^{j_i} \le Q_i}}
    \sum_{\substack{k: \\ B = 2^k \le \frac{Q}{C_1 \ldots C_g} \\ B
        \not\in S}} B (\log 2B)^{A_1} \sum_{\substack{C_i \le c_i \le 2
        C_i \\ (c_i;D)=1}}      \sup_{\partial(\cR) \le M} L(\bc,b(B),\cR).
\end{align*}
Writing $\cR' := \cR/b(B)$, we may now apply our starred level of
distribution formula (Lemma
\ref{thm:lodstar}) to the inner sum, which is bounded from above by
\begin{align*}
&\sum_{\substack{C_i \le c_i \le 2 C_i \\ (c_i;D)=1}}
\sup_{\partial(\cR') \le M/B }
   \left|  \#(\lstarc \cap \cR' \cap \Psi_{b(B)}) - 
  \frac{\rstar(\bc)}{(c_1 \ldots c_g D)^2} \vol(\cR')  \right| \\
& \ll \frac{M}{B} \sqrt{C_1 \ldots C_g} ( \log 2^{g+1} C_1 \ldots
C_g)^{\nu_1'} + C_1 \ldots C_g (\log 2^{g+1} C_1 \ldots C_g)^{\nu_2'}.
\end{align*}
Applying the same reasoning as in \cite{Mar06}, we may then deduce the
level of distribution formula.


\section{Pairs of Forms with Almost Prime Values} \label{chap:almost}

In a  series of papers,  Diamond, Halberstam, and
Richert (and later Diamond and Halberstam) developed a general
multidimensional sieve, which found application \cite{DH97} in the
representation of almost primes by sets of polynomials.  We shall use
the result of Diamond and Halberstam, together with our new level of
distribution formula, to derive similar results concerning the
representation of almost primes by sets of irreducible binary
quadratic forms.  In our application, we will will want to sift the
multiset
$$
\fA := \{ q_1(x,y) \ldots q_g(x,y) : (x,y) \in \Z^2 \cap X\cRo \cap
\Psi\},
$$
and, taking $D=\prod_{p\le 2g}p \prod_{i<j}\Res(q_i,q_j)
\prod_{i=1}^g \delta_i a_i 
c_i$, where $\delta_i$ is the discriminant of the form $q_i(x,y)
:= a_i x^2 + 2b_i xy + c_i y^2$, we define $\Psi := \{ \x
\in \Z^2 : \x \equiv \bz \pmod{D} \}$,
where $\bz$ is chosen such that $(q_i(\bz);D)=1$, for $i=1,\ldots,g$.

The set $\fA$ is sifted by the set $\fP$ of  primes which do not 
divide $D$, and we define $\fPbar$ to be the complement of $\fP$ in
the set of all primes.
We introduce the function $\omega(\cdot)$ which satisfies
$\omega(1)=1$, $\omega(p) =0$ for $p \in \fPbar$, and we let $Y=X^2
\vol(\cRo)/D^2$, an 
approximation to $|\fA|$.  Define $\fA_d := \{ a \in \fA : d|a\}$.
Our level of distribution formula provides information about the
distribution of the $\fA_d$s.  In fact, defining
$$
R_d := |\fA_d| - \frac{\omega(d)}{d} Y, \quad \text{if } \mu(d) \ne 0,
$$
the level of distribution formula shows that the $R_d$ are small on
average, in a sense to be made precise.

Using the above notation, we have the following theorem of Diamond and
Halberstam \cite{DH97}:

\alpheqn

\begin{theorem} \label{thm:hr}

Suppose there exist real constants $\kappa > 1$, $A_1, A_2 \ge 2$, and
$A_3 \ge 1$
such that

\begin{equation} \label{eqn:omega1}
0 \le \omega(p) < p,
\end{equation}

\begin{equation} \label{eqn:omega2star}
\prod_{z_1 \le p < z} \left( 1 - \frac{\omega(p)}{p} \right)^{-1} \le
\left( \frac{\log z}{\log z_1}\right)^\kappa \left( 1 + \frac{A_1}{\log z_1}
  \right), \quad 2 \le z_1 < z,
\end{equation}

\begin{equation} \label{eqn:rkappaalpha}
\sum_{\substack{d < Y^\alpha / (\log Y)^{A_3} \\ (d; \fPbar)=1}}
\mu^2(d) 4^{\nu(d)} | R_d| \le A_2 \frac{Y}{\log^{\kappa+1} Y},
\end{equation}
for some $\alpha$ with $0 < \alpha \le 1$; that

\begin{equation} \label{eqn:acoprimefpbar}
(a ; \fPbar) = 1 \quad \text{for all } a \in \fA,
\end{equation}
and that
\begin{equation} \label{eqn:mualpha}
|a| \le Y^{\alpha \mu} \quad\text{for some }\mu,\text{and for all }
a \in \fA.
\end{equation}

Then there exists a real constant $\beta_\kappa >2$ such that for any real
numbers $u$ and $v$ satisfying
$$
\alpha^{-1} < u < v, \quad \beta_\kappa < \alpha v,
$$
we have
$$
| \{ P_r : P_r \in \fA \} | \gg Y \prod_{p < Y^{1/v}} \left( 1 -
  \frac{\omega(p)}{p} \right),
$$
where
\reseteqn
\begin{equation}
  \label{eqn:rlowerbound}
r > \alpha \mu u - 1 + \frac{\kappa}{f_\kappa(\alpha v)} \int_1^{v/u}
F_\kappa (\alpha v -s) \left(1 - \frac{u}{v}s \right) \frac{ds}{s}.
\end{equation}
The functions $f_\kappa$ and $F_\kappa$ are the solutions to
delay--differential equations specified in \cite{DH97}.  The
parameters $\beta_\kappa$ and $\alpha_\kappa$ that appear in the
delay--differential equations are tabulated in \cite{DHR88} for $1 \le
\kappa \le 10$.
\end{theorem}
\renewcommand{\theequation}{\arabic{equation}}

To apply our level of distribution formula, we must relate $|\fA_{d}|$ to
quantities of the form $\#(\Lambda_{\bc} \cap X \cRo \cap
\Psi)$. Indeed, we have:

\begin{lemma} \label{lem:faformula}
If $d$ is squarefree, then
$$
|\fA_d| = \sum_{\substack{ d | c_1 \ldots c_g \\ c_i |d \\
    i=1,\ldots,g}} \mu(d) \mu(c_1)\ldots \mu(c_g) \# (\Lambda_{\bc}
\cap X \cRo \cap \Psi).
$$
\end{lemma}
For the duration of this proof, let us write $\Omega$ for $\Z^2 \cap X
\cRo \cap \Psi$, then
\begin{align}
|\fA_d| = \# \{\x \in \Omega\} - \#\{\x \in \Omega : d\nmid q_1(\x)
\ldots q_g(\x) \}. \label{eqn:cad1}
\end{align}
Using the fact that $d$ is squarefree, we have
\begin{align*}
& \# \{\x \in \Omega : d\nmid q_1(\x) \ldots q_g(\x) \} = \#
\bigcup_{p|d} \{ \x \in \Omega : (p;q_i(\x))=1, \; i=1,\ldots,g\}  \\
&=-\sum_{\substack{e|d \\ e \ne 1}} \mu(e) \# \{ \x \in \Omega : 
(e;q_i(\x))=1, \; i=1,\ldots,g \},
\end{align*}
where we use the inclusion--exclusion principle in the last line.
Combining this with (\ref{eqn:cad1}), we have
\begin{align}
   \label{eqn:cad2}
|\fA_d| &= \sum_{e|d} \mu(e) \# \{ \x \in \Omega : (e;q_i(\x))=1, \;
i=1,\ldots,g \} 
 \\
&= \sum_{e|d} \mu(e) \sum_{\substack{c_i|e \\
    i=1,\ldots,g}} \mu(c_1)\ldots \mu(c_g) \# (\Lambda_{\bc} \cap
\Omega) \nonumber \\
& = \sum_{\substack{c_i|d \\ i=1,\ldots,g}} \mu(c_1) \ldots \mu(c_g)
\# (\Lambda_{\bc} \cap \Omega) \sum_{\substack{e|d \\ c_i|e
    \\i=1,\ldots,g}} \mu(e). \nonumber
\end{align}

Writing $h := \lcm(c_1,\ldots,c_g)$, and noting that $d$ is
squarefree, we have:
$$ \sum_{\substack{e|d \\ c_i|e \\i=1,\ldots,g}} \mu(e) =
\sum_{\substack{e|d \\ h|e}} \mu(e)= \sum_{\substack{f|d/h}} \mu(fh) =
\sum_{fh|d} \mu(f) \mu(h),$$
This expression equals $\mu(d)$ if $d=h$ and equals
$0$ otherwise.  Now if $d = h$, then $d|c_1 \ldots
c_g$.  Conversely, suppose $d|c_1 \ldots c_g$. As $d$ is a multiple
of $c_1,\ldots,c_g$, we deduce that $h$ divides $d$.  On the
other hand, $d|c_1 \ldots c_g$ and $d$ is squarefree, so $d|
h$, whence $d =h$.  Thus, we may rewrite
equation (\ref{eqn:cad2}) as:
$$
|\fA_d| = \sum_{\substack{c_i|d \\ i=1,\ldots,g \\ d|c_1 \ldots c_g}}
\mu(d) \mu(c_1) \ldots \mu(c_g)
\# (\Lambda_{\bc} \cap \Omega),
$$
as required.

We turn to the problem of defining $\omega(p)$. Recall that for $p \in
\fP$, we would like $Y \omega(p)/p$ to be roughly $|\fA_p|$. Now, we
have the approximation:
$$
\#( \ld \cap \cR \cap \Psi) \approx
\frac{\rho(d_1,\ldots,d_g)}{(d_1 \ldots d_g D)^2}
\vol(\cR).
$$
Taking $\cR := X \cRo$ and bearing in mind Lemma
\ref{lem:faformula}, we choose to define $\omega(p)$ by
\begin{equation} \label{eqn:omegadef}
\omega(p) := p \sum_{\substack{c_i|p \\i=1,\ldots,g \\ p|c_1 \ldots
    c_g}} \mu(p) \mu(c_1)\ldots \mu(c_g)
\frac{\rho(c_1,\ldots,c_g)}{(c_1 \ldots c_g)^2},
\end{equation}
if $p \in \fP$, and $\omega(p)=0$ otherwise.

We'll find it convenient to rewrite $\omega(p)$ as follows:
\begin{lemma} \label{lem:omegaformula}
  For $p \in \fP$,
$$
\omega(p) = p^{-1}( \rho(p,1,\ldots,1)+\ldots+\rho(1,\ldots,1,p)+1-g).
$$
\end{lemma}
Indeed, by equation (\ref{eqn:omegadef}), we have
\begin{align*}
\frac{\omega(p)}{p} &= 1 - \sum_{c_i|p} \mu(c_1)\ldots \mu(c_g)
\frac{\rho(c_1,\ldots,c_g)}{(c_1\ldots c_g)^2} \\
& = 1 - \sum_{0 \le \alpha_i \le 1} (-1)^{\alpha_1+\ldots+\alpha_g}
\frac{\rho(p^{\alpha_1},\ldots,p^{\alpha_g})}{p^{2(\alpha_1+\ldots+\alpha_g)}}.
\end{align*}
By a combination of Lemma \ref{lem:rstarbound} and Lemma
\ref{lem:rbound}, we have that if $(p;\fPbar)=1$, then
$$
\rho(p^{\alpha_1},\ldots,p^{\alpha_g})=p^{2(\alpha_1+\ldots+\alpha_g-1)},
$$
if $0 \le \alpha_i \le 1$ for all $i=1,\ldots,g$ and if $\alpha_i \ne 0$ for at
least two values of $i$.  This leads us to deduce that $\omega(p)/p$ equals
$$
1 - 1 + p^{-2} (
\rho(p,1,\ldots,1)+\ldots+\rho(1,\ldots,1,p)) - p^{-2}\sum_{\substack{0 \le
    \alpha_i \le 1 \\ \exists i \ne j : \alpha_i \alpha_j \ne 0}}
(-1)^{\alpha_1 + \ldots + \alpha_g}.
$$
We evaluate the last sum by taking the sum over all $0 \le \alpha_i
\le 1$, then subtracting the sum over those $\alpha_i$ such that less
that two of the $\alpha_i$ are non-zero.  Thus $\omega(p)/p$ is
equal to:
$$
p^{-2} (\rho(p,1,\ldots,1)+\ldots+\rho(1,\ldots,1,p)) - p^{-2} \sum_{0
\le \alpha_i \le 1}(-1)^{\alpha_1+\ldots+\alpha_g} + p^{-2}(1-g).
$$
The sum over $0 \le \alpha_i \le 1$ is zero, giving the desired
answer.

It remains to verify the conditions of  Theorem (\ref{thm:hr}).
\subsection{Condition (\ref{eqn:omega1})}

In light of Lemma \ref{lem:omegaformula}, we may quickly verify condition
(\ref{eqn:omega1}).  First, we must check that $\omega(p) \ge 0$.
This follows as $\rho(\bd)\ge 1$,for any choice of $\bd$.

On the other hand, in the consideration of the one-form problem, we
proved \cite{Mar06} that if $\alpha_i = \delta_{in}$, then, as
$(p;D)=1$,
$$
\rho(p^{\alpha_1},\ldots,p^{\alpha_g}) = 1 + (p-1)\left(1 +
  \jacobi{\delta_n}{p} \right),
$$
where $\delta_n$ is the discriminant of the quadratic form $q_n$.
So, writing $\chi_n(p) := \jacobi{\delta_n}{p}$,
$$
\omega(p) = g + \sum_{n=1}^g \chi_n(p) - \frac{1}{p}\left( g -1 +
  \sum_{n=1}^g \chi_n(p) \right),
$$
whence $\omega(p) <2g +1/p$.  So $\omega(p)<p$, as we've assumed $p
> 2g$.
\begin{explain}
Incidentally, this inequality explains the factor of $
\prod_{p\le 2g}p$ in our choice of $D$.
\end{explain}

\subsection{Condition (\ref{eqn:omega2star})}

\begin{explain}
This condition expresses the $\kappa$-dimensionality of the sieve
problem.  One should think of the quantity $\omega(p)/p$ as being the
probability that an element of $\fA$ is divisible by $p$, and that
$\kappa$ is the `average' value of $\omega(p)$, in some sense.  In
many sieve problems, one finds that $\kappa=1$, a linear
sieve. However, in our problem, we will demonstrate that $\kappa=g$,
as one would expect from the above definition of $\omega(p)$.

We must prove
$$
\prod_{z_1 \le p < z} \left( 1 - \frac{\omega(p)}{p} \right)^{-1} \le
\left( \frac{\log z}{\log z_1}\right)^\kappa \left( 1 + \frac{A_1}{\log z_1}
  \right), \quad 2 \le z_1 < z.
$$
\end{explain}

Following the argument of \cite{Mar06},
\begin{explain}
we may assume that $z_1 > 2g$, as
$\omega(p)=0$ if $p \le 2g$.  So upon taking logs,
\end{explain}
we must demonstrate
$$
\sum_{z_1 \le p < z} \sum_{i=1}^\infty \frac{\omega(p)^i}{i p^i} \le
\kappa \log \log z - \kappa \log \log z_1 + B_1 / \log z_1,
$$
for $z_1 > 2g$.

The main term is $\sum_{z_1 \le p < z} \omega(p) / p$, which
expands to:
\begin{align*}
&\sum_{z_1 \le p < z} \frac{1}{p}\left(g + \sum_{n=1}^g \chi_n(p)
\right) - \frac{1}{p^2} \left( g-1 + \sum_{n=1}^g \chi_n(p) \right) \\
& = g \log \log z - g \log \log z_1 + \sum_{z_1 \le p < z} \sum_{n=1}^g
\frac{\chi_n(p)}{p} + O (1/z_1) \\
& = g \log \log z - g \log \log z_1 + O (1/ \log z_1).
\end{align*}

\begin{explain}
In estimating the sums involving characters, we use a result of
Mertens', to be found in Chapter 7 of \cite{Dav80}, that for any
non-principal character $\chi$, one has $\sum_{p}
p^{-1}\chi(p)\log p =O(1).$  So
$$
\sum_{z_1 \le p \le z} \frac{\chi(p)}{p} = \sum_{z_1 \le p \le z}
\frac{\chi(p) \log p}{p} \frac{1}{\log p}
\le \sum_{z_1 \le p \le z} \frac{\chi(p) \log p}{p} \frac{1}{\log z_1}
\ll \frac{1}{\log z_1}.
$$
\end{explain}

The error term $\sum_{i=2}^\infty \sum_{z_1 \le p < z}
\omega(p)^i / (i p^i)$ has an upper bound of order $1/ \log z_1$,
as can be seen by using the inequality $\omega(p) \le 2g$ for all
primes $p$.
\begin{explain}
Indeed
\begin{align*}
& \sum_{i=2}^\infty \sum_{z_1 \le p < z} \frac{\omega(p)^i}{i p^i} \le
\sum_{i=2}^\infty \sum_{n \ge z_1} \frac{(2g)^i}{i n^i} \\
& \le \sum_{i=2}^\infty (2g)^i / i \left(\int_{x=z_1}^\infty \frac{1}{x^i} dx +
  \frac{1}{z_1} \right) \ll 1/ z_1 \ll 1/\log z_1,
\end{align*}
as required.
\end{explain}

This completes our verification of condition (\ref{eqn:omega2star}).
We see that $\kappa$, the dimension of the sieve, has the value $\kappa=g$.

\subsection{Condition (\ref{eqn:rkappaalpha})}
\label{subsection:sievelod}

To satisfy this condition, we need a good upper bound for the `average
value' of
$$
|R_d| := \left| |\fA_d| - \frac{\omega(d)}{d} Y \right|,
$$
for squarefree $d$.
\begin{explain}
Essentially, we shall sum $|R_d|$ as $d$ varies
in some range.  In this problem, the range of summation is referred to
as the level of distribution, and it is our aim is to ensure that the
level of distribution is as large as possible, whilst requiring that
the sum be bounded above by $Y / (\log Y)^{g+1}$.
\end{explain}%
We've already derived a formula for $|\fA_d|$ in Lemma
\ref{lem:faformula}, so we seek an equivalent formula for
$\omega(p)$.  We make use of the following auxiliary Lemma:

\begin{lemma}  \label{lem:multlem}
  Let $t$ be a function of $g$ variables such that
$$t(d_1,\ldots,d_g)t(e_1,\ldots,e_g)=t(d_1 e_1,\ldots, d_g e_g)$$
if $d_i$ and $e_i$ are squarefree numbers and if $(d_i;e_i)=1$ for
$i=1,\ldots,g$. If the function $f$ is defined by
$$
f(p) := \sum_{\substack{c_i | p \\ p | c_1 \ldots c_g}}
t(c_1,\ldots,c_g),
$$
for every prime $p$ in a set $S$, then
$$
\prod_{p|d}f(p) = \sum_{\substack{c_i|d \\ d | c_1 \ldots
    c_g}} t(c_1,\ldots,c_g),
$$
if $d$ is a squarefree number, all of whose prime factors lie in $S$.
\end{lemma}
The idea of the proof is to define
$$
T_n := \{ (c_1,\ldots,c_g) \in \N^g : n | c_1 \ldots c_g,\text{ and } c_i | n,
\text{ for } i=1,\ldots,g \}.
$$ 
The mapping from $\prod_{p|d} T_p$ to $T_d$ defined by
$\prod_{p|d} (c_{1p},\ldots, c_{gp}) \mapsto (c_1,\ldots,c_g)$, where
$c_i := \prod_{p|d} c_{ip}$ is a bijection.  The Lemma then follows
easily.

Applying Lemma \ref{lem:multlem} to the definition of $\omega(p)$, we
deduce:
$$
\frac{\omega(d)}{d} = \prod_{p|d} \frac{\omega(p)}{p} =
\sum_{\substack{c_i | d \\ d | c_1 \ldots c_g}} \mu(d) \mu(c_1)\ldots
\mu(c_g) \frac{\rho(c_1,\ldots,c_g)}{(c_1 \ldots c_g)^2},
$$
whence
$$
|R_d| \le \sum_{\substack{c_i | d \\ d | c_1 \ldots c_g}} \left|
  \#(\Lambda_{\bc} \cap X \cRo \cap
    \Psi) - Y \frac{\rho(\bc)}{(c_1 \ldots c_g)^2} \right|.
$$

As we consider the sum in condition (\ref{eqn:rkappaalpha}), there is
some flexibility in the choice of the constant $\alpha$.
Indeed, we aim to show that for any $\alpha<1$, the sum in question is
bounded from above by $Y / (\log Y)^{g+1}$.  To this end, we take
$\alpha<1$ and consider the sum:
$$
E := \sum_{\substack{d < Y^\alpha \\ (d;\fPbar)=1}} \mu^2(d)
    4^{\nu(d)} \sum_{\substack{c_i | d \\ d | c_1 \ldots c_g}} \left|
      \#(\Lambda_{\bc} \cap X \cRo \cap
    \Psi) - Y \frac{\rho(\bc)}{(c_1 \ldots c_g)^2} \right|.
$$
Let us introduce another
variable, $k$, which specifies the least common multiple of the
$k_{ij}$, where $k_{ij}:= (c_i;c_j)$.  Define the sets $U$ and $T_k$ by
$$
U := \{ d \in \Z : (d;\fPbar)=1, \mu^2(d)=1\},
$$
and
$$
T_k := \{ (c_1, \ldots,c_g) \in U^g : k = \lcm_{i<j} k_{ij} \}.
$$
As $d$ is squarefree, the conditions $d|c_1 \ldots c_g$ and $c_i|d$
for $i=1,\ldots,g$ are equivalent to $d=\lcm(c_1,\ldots,c_g)$.
Moreover, $\lcm_{i<j} k_{ij}$ divides
$\lcm(c_1,\ldots,c_g)$, so $k|d$, hence $k<Y^\alpha$. This
leads to the expression
$$
E =  \sum_{\substack{d < Y^\alpha \\ d \in U }} \mu^2(d)
4^{\nu(d)} \sum_{\substack{k < Y^\alpha \\ k \in U}}
\sum_{\substack{(c_1, \ldots, c_g) \in T_k \\
  \lcm(c_1,\ldots,c_g)=d }}
\left| \#(\Lambda_{\bc} \cap X \cRo \cap
    \Psi) - Y \frac{\rho(\bc)}{(c_1 \ldots c_g)^2} \right|.
$$

Our guiding principle is that the main term in $E$ arises from those
$c_1,\ldots,c_g$ which are mutually pairwise coprime.  We make the
following change of variable:
$$
f_i := \frac{c_i}{(c_i;k)} \quad \text{ for } i=1,\ldots,g.
$$
Then the $f_1,\ldots,f_g$ are mutually pairwise coprime.  \begin{explain}
Indeed if $p$ is a prime divisor of $f_i$ and $f_j$ then $p|c_i,c_j$.
However, $k$ is
squarefree, so $p\nmid k$, so $p\nmid (c_i;c_j)$, a contradiction.
\end{explain}

How does the change of variable affect the summand?  We'll show that
the map $ \x \to k^{-1} \x$ is a bijection from $\Lambda_{\bc} \cap X
\cRo \cap \Psi$ to $\Lambda_{\bdf} \cap k^{-1}X \cRo \cap \Psi_k$.  It
suffices to show that the map is a bijection from $\Lambda_{\bc}$ to
$\Lambda_{\bdf}$.

First, we demonstrate that if $\x \in
\Lambda_{\bc}$, then $\x \equiv 0 \smod{k}$.  Indeed, if $\x \in
\Lambda_{\bc}$, then for all $i<j$, one has $c_i | q_i(\x)$ and $c_j |
q_j(\x)$, so we have the simultaneous equations:
$$
q_i(\x) \equiv q_j(\x) \equiv 0 \smod{k_{ij}}.
$$
However, $k_{ij}$ is coprime to the resolvent of $q_i$ and $q_j$, so
we deduce $\x \equiv 0 \smod{k_{ij}}$.  As this is true for all $i<j$,
we deduce that $\x \equiv 0 \smod{k}$.

We may write $\x = k \by$, for some $\by \in \Z^2$, so for all
$i=1,\ldots,g$, we have $c_i | k^2 q_i(\by)$, and $f_i|c_i$, so
$f_i|k^2 q_i(\by)$, but $(k;f_i)=1$, so $f_i|q_i(\by)$.  Thus $\by \in
\Lambda_{\bdf}$.  Conversely, if $\by \in \Lambda_{\bdf}$, one easily
checks that $k \by \in \Lambda_{\bc}$, completing the proof of
bijectivity.

In a similar manner, we will treat the quantity $\rho$. For all $i$,
we have $c_i = (c_i;k)f_i$ is squarefree, so $(c_i;k)$ and $f_i$ are
coprime;  by multiplicativity of $\rho$, we find:
$$
\frac{\rho(c_1,\ldots,c_g)}  {(c_1 \ldots c_g)^2} =
\frac{\rho((c_1;k),\ldots,(c_g;k))}{(c_1;k)^2 \ldots (c_g;k)^2}
\frac{\rho(f_1,\ldots, f_g)}{(f_1\ldots f_g)^2}.
$$

Recall that $\rho((c_1;k),\ldots,(c_g;k))$ counts those $\x$ in
$[0,(c_1;k)\ldots (c_g;k)^2)$ such that $q_i(\x) \equiv 0
\smod{(c_i;k)}$, for $i=1,\ldots,g$.

Suppose that $i<j$. Note that $(c_i;c_j;k)=((c_i;c_j);k)=(c_i;c_j)=k_{ij}$,
so if $(c_i;k) | q_i(\x)$ and
$(c_j;k)|q_j(\x)$, then $q_i(\x) \equiv q_j(\x) \equiv 0
\smod{k_{ij}}$.  As before, we have $\x \equiv 0 \smod{k_{ij}}$, and
  indeed $\x \equiv 0 \smod{k}$.  Conversely, if $\x \equiv 0
  \smod{k}$, then for $i=1,\ldots,g$, we have $q_i(\x) \equiv 0
  \smod{k}$; therefore, $(c_i;k) | q_i(\x)$.

We deduce that $\rho((c_1;k),\ldots,(c_g;k)) = k^{-2} (c_1;k)^2\ldots
(c_g;k)^2$, and hence
$$
\frac{\rho(c_1,\ldots,c_g)}{(c_1\ldots c_g)^2} =
\frac{\rho(f_1,\ldots,f_g)}{(k f_1\ldots f_g)^2}.
$$

As a result of these formul\ae, we may write
$
\left| \# (\Lambda_{\bc} \cap X \cRo \cap \Psi) - Y
  \frac{\rho(\bc)}{(c_1\ldots c_g)^2}\right|$ as $\left| \# (
  \Lambda_{\bdf} \cap k^{-1} X \cRo \cap \Psi_k) - Y
  \frac{\rho(\bdf)} {(k f_1\ldots f_g)^2} \right|.
$
It remains to reformulate the condition $d= \lcm(c_1,\ldots, c_g)$.
Now $\lcm(c_1,\ldots,c_g) = \lcm_{1\le i \le g}(c_i;k)f_i$.  We
have already observed that $(c_i;k)$ is coprime to $f_i$ for all $i$,
and that $f_1,\ldots,f_g$ are mutually pairwise coprime, so
$\lcm(c_1,\ldots,c_g)= f_1\ldots f_g \lcm_{1\le i \le g}(c_i;k).$
One has $k = \lcm_{1 \le i \le g}(c_i;k)$.
\begin{explain}
  Indeed, if $p$ is a prime divisor of $k$ then $p$ divides
  $(c_i;c_j)$ for some $i<j$, thus $p | (c_i;k)$.  Hence $p | \lcm_{1
    \le i \le g} (c_i;k)$, and we deduce that $k| \lcm_{1 \le i \le
    g}(c_i;k)$ .  Conversely, if $p|\lcm_{1 \le i \le
    g}(c_i;k)$, then $p | (c_i;k)$ for some $i$. In particular,
  $p|k$.  We deduce that $\lcm_{1\le i \le g}(c_i;k)$ divides $k$, as required
\end{explain}

So we may write the error term as
\begin{align*}
E & \le \sum_{\substack{d < Y^\alpha \\ d \in U}} \mu^2(d) 4^{\nu(d)}
\sum_{\substack{k < Y^\alpha \\ k \in U}}
\sum_{\substack{(f_1,\ldots,f_g) \in T_1 \\ k f_1 \ldots f_g = d}}
\left| \#(\Lambda_{\bdf} \cap k^{-1} X \cRo \cap \Psi_k) - Y
  \frac{\rho(\bdf)}{(k f_1\ldots f_g)^2} \right| \\
& \le \sum_{\substack{k<Y^\alpha \\ k \in U}}
\sum_{\substack{(f_1,\ldots,f_g)\in T_1 \\ kf_1 \ldots f_g <
    Y^\alpha}} \left| \#(\Lambda_{\bdf} \cap k^{-1} X \cRo \cap \Psi_k) - Y
  \frac{\rho(\bdf)}{(k f_1\ldots f_g)^2} \right|  4^{\nu(k f_1 \ldots f_g)}.
\end{align*}
Now for all $\eps >0$, and for all $m \in \N$, one has $4^{\nu(m)}
\ll_\eps m^\eps$, so $4^{\nu(k f_1 \ldots f_g)} \ll_\eps Y^\eps$, as
$\alpha\le 1$.  Thus,
$$
E \ll_\eps Y^\eps \sum_{\substack{k<Y^\alpha \\ k \in U}}
\sum_{\substack{(f_1,\ldots,f_g)\in T_1 \\ kf_1 \ldots f_g <
    Y^\alpha}} \left| \#(\Lambda_{\bdf} \cap k^{-1} X \cRo \cap \Psi_k) - Y
  \frac{\rho(\bdf)}{(k f_1\ldots f_g)^2} \right|.
$$

Let $E(k)$ denote the inner sum.  We evaluate $E(k)$ by splitting the
range for $f_1,\ldots,f_g$ into dyadic intervals.  For $i=1,\ldots,g$,
there exists integers $n_1,\ldots,n_g$ such that $2^{n_i-1} \le f_i <
2^{n_i}$, so $2^{n_1+\ldots+n_g-g} \le f_1 \ldots f_g < Y^\alpha/k$.
Thus
$$
E(k) \le \sum_{\substack{n_1,\ldots, n_g: \\
    k 2^{n_1+\ldots+n_g-g}<Y^\alpha}} \sum_{\substack{f_i <
    2^{n_i}\\i=1,\ldots,g}} \left| \#(\Lambda_{\bdf} \cap k^{-1} X
  \cRo \cap \Psi_k) - Y \frac{\rho(\bdf)}{(k f_1\ldots f_g)^2} \right|.
$$

Applying the level of distribution formula to the inner sum, $E(k)$ is
bounded above by a quantity of order
$$
\sum_{\substack{n_1,\ldots, n_g: \\ k 2^{n_1+\ldots+n_g-g} <
    Y^\alpha}} 2^{\sum_{i} n_i} (\log 2^{\sum_i n_i +1})^{\nu_1} + 
\frac{Y^{1/2}}{k}(2^{\sum_{i}n_i})^{1/2}(\log 2^{\sum_i n_i + 1})^{\nu_2}.
$$
To calculate this, we introduce the quantity $Q := \log_2 (2^g Y^\alpha /
k)$.  Then the condition $0\le n_1+\ldots+n_g<Q$ is equivalent to the
conditions $0\le n_i <Q$, for $1 \le i \le g-1$, and  $0 \le n_g < Q -
n_1-\ldots-n_{g-1}$.  Hence we split the
sum for $E(k)$ to give an upper bound of
$$
\sum_{n_1 < Q} \ldots \sum_{n_g <
  Q-n_1-\ldots-n_{g-1}} 2^{\sum_i n_i}
(1+\sum_i n_i)^{\nu_1} + \frac{Y^{1/2}}{k}(2^{\sum_i
  n_i})^{1/2}(1+\sum_i n_i)^{\nu_2}.
$$
Now if $n_i$ are nonnegative, then $(1+\sum_{i=1}^g n_i) \le 
\prod_{i=1}^g (1+n_i)$.
\begin{explain}
  The result is trivial if $g=1$.  Suppose then that
  $(1+\sum_{i=1}^{g-1} n_i) \le \prod_{i=1}^{g-1} (1+n_i)$.  Then
$$
(1+\sum_{i=1}^g n_i) \le n_g + \prod_{i=1}^{g-1}(1+n_i) \le
n_g\prod_{i=1}^{g-1} (1+n_i) + \prod_{i=1}^{g-1}(1+n_i) = \prod_{i=1}^g
(1+n_i),
$$
as required.
\end{explain}
Using this bound gives the following upper bound for $E(k)$:
\begin{align*}
&\sum_{n_1 < Q} 2^{n_1}(1+n_1)^{\nu_1} \ldots \sum_{n_g <
  Q-n_1-\ldots-n_{g-1}} 2^{n_g} (1+n_g)^{\nu_1} \\ &+\frac{Y^{1/2}}{k}
\sum_{n_1<Q} 2^{n_1/2}(1+n_1)^{\nu_2} \ldots \sum_{n_g <
  Q-n_1-\ldots-n_{g-1}} 2^{n_g/2}(1+n_g)^{\nu_2}.
\end{align*}

Following a similar line to \cite{Mar06}, we estimate $E(k)$:
\begin{explain}
Now if $\beta > 0$, $\theta \ge 1$, and $N \ge 1$, then
$$
\sum_{0 \le t < N} 2^{t \beta} (t+1)^\theta \ll_\beta N^\theta
2^{ N \beta}.
$$
Applying this result to our estimate for $E(k)$, we arrive at:
$$
E(k)  \ll 2^Q Q^{\nu'}
+ \frac{(2^Q Y)^{1/2}}{k} Q^{\nu'},
$$
for some absolute constant $\nu'$.
Recalling the definition of $Q$, we have the upper bound:
\end{explain}
$$
E(k) \ll \frac{Y^\alpha}{k} (\log Y)^{\nu'} +
\frac{Y^{(\alpha+1)/2}}{k^{3/2}} (\log Y)^{\nu'}.
$$
Finally, we sum $E(k)$ over k:
$$
E \ll_\eps Y^\eps Y^\alpha (\log Y)^{\nu'+1} + Y^\eps Y^{(\alpha+1)/2}
(\log Y)^{\nu'} \ll_\eps Y^{\max(\alpha+\eps,\alpha/2+1/2+\eps)}.
$$
If we choose $\eps = (1-\alpha)/3$, then $E \ll_\alpha Y / (\log
Y)^{\kappa+1}$, and condition (\ref{eqn:rkappaalpha}) is
satisfied for any $\alpha <1$.

\subsection{Conditions (\ref{eqn:acoprimefpbar}) and (\ref{eqn:mualpha})}

Every $a \in \fA$ can be represented as $a=q_1(\x) \ldots q_g(\x)$,
with $\x \in \Psi$. For all $\x \in
\Psi$, and
for $i=1,\ldots,g$, we have
$(q_i(\x);D)=1$, so $(a;D)=1$, whence $(a;\fPbar)=1$,
satisfying condition (\ref{eqn:acoprimefpbar}).

Now observe
that for all $a \in \fA$, one has $|a| \ll X^{2g} \ll Y^g$.  That is,
there exists a constant $C$ (depending on the choice of forms
$q_i$ and the region $\cRo$) such that $|a| \le C Y^g$ for all $a \in
\fA$.  Define $\theta$ by $C=Y^\theta$.  In order to satisfy
condition (\ref{eqn:mualpha}), we need $|a| \le Y^{\mu \alpha}$, and
it is sufficient to choose $\alpha <1$ and $\mu$ such that $\mu \ge
(g+\theta)/\alpha$.

The lower bound for $r$ given in (\ref{eqn:rlowerbound}) is continuous in
$\alpha$ and $\mu$, so we may take $\alpha=1$ and $\mu=g$ if we can
find $\alpha$, $\mu$ satisfying condition
(\ref{eqn:mualpha}) such that $|\alpha-1|,
|\mu-g|<\eta$ for any $\eta>0$.  Indeed, set
$\mu = (g+\theta)/\alpha$.  For $\alpha<1$, the above condition
translates into $\alpha > (g+\theta)/(g+\eta)$ and $\alpha > 1-\eta$.
We can choose such a value of $\alpha$ provided that $\theta<\eta$.
Now $\theta=\log C / \log Y$, so the
condition will be satisfied for all sufficiently large $Y$.  

\subsection{Application of Theorem \ref{thm:hr} }

Our work up to this point has been in order to justify taking the
parameters $\alpha=1$ and $\kappa=\mu=g$ in the lower bound
(\ref{eqn:rlowerbound}).  To compute an optimal value for the lower
bound, we solve the delay--differential equations of \cite{DH97} and
minimise the result over $u$ and $v$.  This is all carried out using
Mathematica \cite{Wol88}.

We deduce that
for sufficiently large $X$, there exists a constant $v_\kappa > 2$ such
that
$$
|\{P_{r_M} : P_{r_M} \in \fA \} | \gg X^2 \prod_{p < X^{2/v_\kappa}} \left( 1 -
  \frac{\omega(p)}{p} \right),
$$
with $r_M$ as in Table \ref{table:r_M}. This proves Theorem \ref{thm:main}.

\begin{explain}
\begin{tabular}{r r r r r}
$\kappa$ & \multicolumn{1}{c}{$\alpha_\kappa$} &
\multicolumn{1}{c}{$\beta_\kappa$} &
$r_{DH}$ & $r_M$
\\ 
\hline
2 &  5.3577.. &  4.2644.. & 7  &  5 \\
3 &  8.3719.. &  6.6408.. & 12 &  8 \\
4 & 11.5317.. &  9.0722.. & 17 & 12 \\
5 & 14.7735.. & 11.5347.. & 23 & 16 \\
6 & 18.0679.. & 14.0146.. & 29 & 20  \\
7 & 21.3989.. & 16.5042.. & 35 & 25  \\
8 & 24.7571.. & 18.9988.. & 41  & 29  \\
9 & 28.1326.. & 21.4955.. & 47  & 34 \\
10 & 31.5320.. & 23.9924..& 54  & 39
\end{tabular}
\end{explain}


\begin{acknowledgements}
I would like to thank Roger Heath-Brown for many useful conversations,
and for guiding me to the right questions.  
\end{acknowledgements}

\bibliography{bible}
\bibliographystyle{amsplain}

\begin{flushleft}
Mathematical Institute, \\
University of Oxford, \\
24--29 St Giles', \\
Oxford, \\
OX1 3LB \\
Email: {\tt marasing@maths.ox.ac.uk} \\
\end{flushleft}
\end{document}